\documentclass[openany]{book}
\usepackage[german]{babel}
\usepackage[all,cmtip]{xy}
\usepackage{hyperref}
\usepackage{amsmath,amssymb,array,fancyhdr,fp,graphicx,
  ifthen,latexsym,multicol,pifont,relsize,stmaryrd,textcomp,yfonts}
\begin{document}
\newcommand {\AAM} {\mathbb{A}}
\newcommand {\BB} {\mathbb{B}}
\newcommand {\CC} {\mathbb{C}}
\newcommand {\DD} {\mathbb{D}}
\newcommand {\EE} {\mathbb{E}}
\newcommand {\FF} {\mathbb{F}}
\newcommand {\HH} {\mathbb{H}}
\newcommand {\II} {\mathbb{I}}
\newcommand {\KK} {\mathbb{K}}
\newcommand {\MM} {\mathbb{M}}
\newcommand {\NN} {\mathbb{N}}
\newcommand {\PP} {\mathbb{P}}
\newcommand {\QQ} {\mathbb{Q}}
\newcommand {\RR} {\mathbb{R}}
\newcommand {\TT} {\mathbb{T}}
\newcommand {\YY} {\mathbb{Y}}
\newcommand {\ZZ} {\mathbb{Z}}
\newcommand {\AAA} {\mathcal{A}}
\newcommand {\BBB} {\mathcal{B}}
\newcommand {\CCC} {\mathcal{C}}
\newcommand {\DDD} {\mathcal{D}}
\newcommand {\EEE} {\mathcal{E}}
\newcommand {\FFF} {\mathcal{F}}
\newcommand {\GGG} {\mathcal{G}}
\newcommand {\HHH} {\mathcal{H}}
\newcommand {\III} {\mathcal{I}}
\newcommand {\JJJ} {\mathcal{J}}
\newcommand {\KKK} {\mathcal{K}}
\newcommand {\LLL} {\mathcal{L}}
\newcommand {\MMM} {\mathcal{M}}
\newcommand {\NNN} {\mathcal{N}}
\newcommand {\OOO} {\mathcal{O}}
\newcommand {\PPP} {\mathcal{P}}
\newcommand {\RRR} {\mathcal{R}}
\newcommand {\SSS} {\mathcal{S}}
\newcommand {\TTT} {\mathcal{T}}
\newcommand {\UUU} {\mathcal{U}}
\newcommand {\VVV} {\mathcal{V}}
\newcommand {\XXX} {\mathcal{X}}
\newcommand {\ZZZ} {\mathcal{Z}}
\newcommand {\Ab}[2][1400]
{\subsubsection[{\Century{#1}#2}]
{\NIF{\Century{#1}#2}}}

\newcommand {\Abo}[1] {\bigskip\NIFG{#1}\m}

\newenvironment {Abstract} [1][Abstract]
{\begin{center}{\Century{1000}\F{#1}}\\[1em]
\begin{minipage}{9cm}\Century{900}} {\end{minipage}\end{center}}

\newcommand {\aitem} [3] {\item [#1] \hypertarget{#2}{} \F{#3:}}

\newcommand {\Align}[1]
{\begin{displaymath}\begin{aligned}#1\end{aligned}\end{displaymath}}

\newenvironment {biblio} {\begin{list}{}
{\renewcommand{\makelabel}[1]{\hfill[##1]}}} {\end{list}}

\newcommand {\CK}[1] {{\textsmaller[2]{#1}}}
\newcommand {\ein}[1] {\hspace*{#1em}} 
\newcommand {\F}[1] {\textbf{#1}}
\newcommand {\Filt} {\operatorname {Filt}}
\newcommand {\Frac}[2] {\displaystyle\frac{#1}{#2}}

\newcommand {\gauto} [3][1.5]
{{\renewcommand{\arraystretch}{#1}
\begin{tabular} [t] {|>{$}l<{$}|p{#2cm}|}
\hline&\\[-1.5em] #3 \\[0.5em]\hline\end{tabular}}}

\newcommand {\gd} {\;\Longleftrightarrow\;}
\newcommand {\gk}[1] {\{#1\}}
\newcommand {\hme} {^{-1}}
\newcommand {\K}[1] {{\itshape #1}}

\newcommand {\kapitel} [1] {\setcounter{Kapzaehler}{#1}\setcounter{Satzzaehler}{1}}

\newenvironment {keywords} [1][Keywords]
{\Century{950}\NI\K{Keywords: }} {}

\newcommand {\leer} {=\emptyset}
\newcommand {\M} [1] {{\ensuremath{#1}}}
\newcommand {\m} {\medskip}
\newcommand {\Mal} {\lim\limits}
\newcommand {\Maln}[1][x] {\Mal_{#1 \Mpk 0}}
\newcommand {\Malu}[1][x] {\Mal_{#1 \Mpk \infty}}
\newcommand {\Malun} {\Malu[n]}
\newcommand {\Map} {\prod\limits}
\newcommand {\Mas} {\sum\limits}
\newcommand {\Matrix}[1] {\begin{pmatrix}#1\end{pmatrix}}
\newcommand {\Max} {\on{Max}}
\newcommand {\Mf} {\mathop{\bigcirc}\limits}
\newcommand {\Mfa}[5] [-0.6]
{\raisebox{#1em}{\M{\begin{array}{r@{\;}c@{\;}l}#2 & \Mp & #3\\
  #4 & \longmapsto & #5\end{array}}}}

\newcommand {\MfDreifaelle}[6] {\M{\begin{cases}{#1} & {\quad #2}\\
  {#3} & {\quad #4}\\ {#5} & {\quad #6}\end{cases}}}

\newcommand {\MfZweifaelle}[5][0ex] {\M{\begin{cases}{#2} & {\quad #3}\\[#1]
  {#4} & {\quad #5}\end{cases}}}

\newcommand {\Min} {\on{Min}}
\newcommand {\Mk}[3] {\M{\left#1 #3 \right#2}}
\newcommand {\Mkr}[1] {\Mk {(} {)} {#1}}
\newcommand {\mm} {\bigskip}
\newcommand {\Mod} {\operatorname {mod}}
\newcommand {\Modg} {\bigcap\limits}
\newcommand {\Moi} {\operatorname {int}}
\newcommand {\Movd} {\sqcup}
\newcommand {\Movgd} {\bigsqcup\limits}
\newcommand {\Movg} {\bigcup\limits}
\newcommand {\Mp} {\M{\longrightarrow}}
\newcommand {\Mpb} [1] {\M{\overset{#1}{\Mp}}}
\newcommand {\Mpd} {\M{\longleftrightarrow}}
\newcommand {\Mpdk} {\M{\leftrightarrow}}
\newcommand {\Mpi} {\M{\;\Longrightarrow\;}}
\newcommand {\Mpk} {\M{\rightarrow}}
\newcommand {\NI} {\noindent}
\newcommand {\NIF}[1] {\NI\F{#1}}
\newcommand {\NIFG}[1] {\NI\F{\relsize{1}#1}}
\newcommand {\nleer} {\ne\emptyset}

\newcommand {\Nummer} [2][ ]
  {\ifthenelse{\equal{#2}{}}{\Nummernzaehlung[#1].}{\Nummernzaehlung[#1]{} #2.}}

\newcommand {\Nummernzaehlung} [1][ ]
  {\ifthenelse{\value{Kapzaehler}=-1}{}
  {\ifthenelse{\value{Kapzaehler}=0}
  {#1\arabic{Satzzaehler}\stepcounter{Satzzaehler}}
  {#1\arabic{Kapzaehler}.\arabic{Satzzaehler}\stepcounter{Satzzaehler}}}}

\newcommand {\oben}[1] {\raisebox{1ex}{\CK{#1}}}
\newcommand {\unten}[1] {\raisebox{-0.5ex}{\CK{#1}}}
\newcommand {\on}[1] {\operatorname{#1}}
\newcommand {\set} {\setlength}
\newcommand {\sm} {\smallskip}
\newcommand {\U}[1] {\underline{#1}}

\newcommand {\zitat} [4][]
{#3\ifthenelse{\equal{#3}{}}{}{ }
[\hyperlink{#2}{\F{#4}}\ifthenelse{\equal{#1}{}}{}{, #1}]}

\newcommand {\zweim} {\hspace*{2em}}
 \newcounter{Kapzaehler}
\newcounter{Satzzaehler}
\newcommand {\Comment}[2][] {\NI\F{Comment\Nummer{#1}} #2\m}
\newcommand {\Conjecture}[2][] {\NI\F{Conjecture\Nummer{#1}} #2\m}
\newcommand {\Corollary}[2][] {\NI\F{Corollary\Nummer{#1}} \K{#2}\m}
\newcommand {\Corollarynk}[2][] {\NI\F{Corollary\Nummer{#1}} #2\m}
\newcommand {\Definition}[2][] {\NI\F{Definition\Nummer{#1}} #2\m}
\newcommand {\Definitionz}[1][] {\NI\F{Definition\Nummer{#1}}\ }
\newcommand {\Example}[2][] {\NI\F{Example\Nummer{#1}} #2\m}
\newcommand {\Examples}[2][] {\NI\F{Examples\Nummer{#1}} #2\m}
\newcommand {\Examplesz}[1][] {\NI\F{Examples\Nummer{#1}}\ }
\newcommand {\Examplez}[1][] {\NI\F{Example\Nummer{#1}}\ }
\newcommand {\Exercise}[2][] {\NI\F{Exercise\Nummer{#1}} #2\m}
\newcommand {\Exercisez}[1][] {\NI\F{Exercise\Nummer{#1}}\ }
\newcommand {\Lemma}[2][] {\NI\F{Lemma\Nummer{#1}} \K{#2}\m}
\newcommand {\Missing}[1][] {\framebox{MISSING {#1}}}
\newcommand {\Notation}[2][] {\NI\F{Notation\Nummer{#1}} #2\m}
\newcommand {\Note}[2][] {\NI\F{Note\Nummer{#1}} #2\m}
\newcommand {\Notez}[1][] {\NI\F{Note\Nummer{#1}}\ }
\newcommand {\Remark}[2][] {\NI\F{Remark\Nummer{#1}} #2\m}
\newcommand {\Remarkz}[1][] {\NI\F{Remark\Nummer{#1}}\ }
\newcommand {\Proof}[1] {\U{Proof.} #1\m}
\newcommand {\Proposition}[2][] {\NI\F{Proposition\Nummer{#1}} \K{#2}\m}
\newcommand {\Reference}[2][] {\NI\F{Reference\Nummer{#1}} #2\m}
\newcommand {\Standing}[2][] {\NI\F{Standing hypothesis\Nummer{#1}} #2\m}
\newcommand {\Test}[2][] {\NI\F{Test\Nummer{#1}} #2\m}
\newcommand {\Testz}[1][] {\NI\F{Test\Nummer{#1}}\ }
\newcommand {\Theorem}[2][] {\NI\F{Theorem\Nummer{#1}} \K{#2}\m}
\newlength{\HPT} \set{\HPT}{0.01pt}

\newlength{\HAB} \set{\HAB}{0.0118pt}

\newcommand {\defont}[2] {\newcommand {#1}[1]{\fontfamily{#2}%
\fontsize{##1\HPT}{##1\HAB}\selectfont}}

\defont {\Century}{pnc}
\newcommand {\DV}[4]
{\Frac{#1 \mid #2}{#3\mid #4}}

\newcommand {\Reg} {\on{Reg}}

\newcommand {\Sing} {\on{Sing}}

\newcommand {\Top} {\on {Top}}
\newcommand {\abg} {\overset{_\text{---}}{\subset}}

\Century{1100}
\vspace*{2.5cm}
\begin{center}
\F{\Century{1800}The topological resolution of a\\
finite closure space}
\end{center}\mm

\begin{center}
Josef Eschgf\"aller\\[2ex]
{\small Universita' degli Studi di Ferrara (retired)\\[2ex]
esg\symbol{64}unife.it}
\end{center}\mm

\begin{Abstract}
\NI For every finite closure space $X$ one can define
a finite topological space $\Top X$ together with a natural
projection $\Top X\Mp X$. This could allow to apply
the techniques of topological combinatorics to the study of
finite closure spaces.
\end{Abstract}

\kapitel{1}
\Ab{1. Preliminaries}

\NI For a set $X$ we denote by $\PPP(X)$ the power set of $X$ (set of all subsets of $X$) and
by $\PPP_*(X)$ the set of all non-empty subsets of $X$.

Let $\Filt X$ be the set of all filters on $X$ and, for a subset $A\subset X$, define\m

$\widehat{A} := \gk{B\subset X \mid A\subset B}$\m

\NI Then, if $X$ is finite, it is well known and immediate to show that there exists a natural bijection
$\Filt X\Mpd \PPP_*(X)$ which sends a filter to the intersection of its elements and
a non-empty subset $A\subset X$ to $\widehat{A}$.

This motivates Definition 3.25.

For a topological space $X$ and a point $x\in X$ we denote by $\UUU(x)$ the set of all
neighborhoods of $x$. We use the same notation for the neighborhoods in a closure space.

The elements of the topological resolution $\Top X$ are ordered pairs $(x,M)$.\\
In order to shorten the notation, we shall denote such a pair by $xM$.

For a finite quasiordered set $(T,\le)$ and $t\in T$ we denote by\m

$U_t:=\gk{s\in T \mid s\ge t}$\m

\NI the upper set determined by $t$ which is at the same time the smallest neighborhood
of $t$ if we consider $T$ as a topological space (cf. Proposition 2.1).
\kapitel{2}
\Ab{2. Finite topological spaces}

\Proposition {(1) Let $T$ be a topological space. Then we may\\introduce a quasiordering on $T$
by defining\m

$t\le s\gd t\in\overline{s}$\m

(2) Viceversa, if $(T,\le)$ is a quasiordered set, then we obtain a topology on $T$ if we define\m

$\UUU(t) := \widehat{U_t}=\gk{V\subset T \mid V\supset U_t}$\m

\NI where, as in the preliminaries, $U_t:=\gk{s\in T \mid s\ge t}$.

Notice that in this way every point $t$ has a smallest neighborhood which coincides with $U_t$.

It is also immediate that $t\le s\gd s\in U_t \gd U_s\subset U_t$.

(3) If $T$ is finite, the constructions in (1) and (2) are one the reversal of the other, so that the
concepts of finite topological space and of finite quasi\-ordered set coincide.\m

(4) $(T,\UUU)$ is T\unten{0} iff $(T,\le)$ is partially ordered.\m

(5) A mapping between finite topological spaces is continuous iff it is order preserving.}

\Proof {This is well known, see e.g. Birkhoff [1, p. 117], Ern\'e [4], Stong [15],
and (with reversed ordering) Barmak [9, p. 2-3], May [14, p. 3].

For a comprehensive exposition of the algebraic topology of finite topological spaces (and hence of finite
quasiordered sets) see Barmak [9].}
\kapitel{3}
\Ab{3. Finite closure spaces}

\Definition {Let $X$ be a set and $^{-}:\PPP(X)\Mp\PPP(X)$ be a mapping such that
for every $A,B\subset X$ the following conditions are satisfied:\m

(1) $A\subset\overline{A}$.

(2) $A\subset B\Mpi \overline{A}\subset \overline{B}$.

(3) $\overline{\overline{A}}=\overline{A}$.\m

\NI $X=(X,^{-})$ is then called a \K{closure space}.}

\Standing {Let $X,Y,Z,...$ be \K{finite} closure spaces.}

\Definition {A point $x\in X$ is \K{inessential}, if $x\in\overline{\emptyset}$.

Otherwise $x$ is said to be \K{essential}.}

\Definition {A subset $A\subset X$ is \K{closed} if $\overline{A}=A$.}

\Definition {A subset $U\subset X$ is \K{open} if $X\setminus U$ is closed.}

\Remark {A subset $A\subset X$ is closed iff there exists $B\subset X$\\such that
$A=\overline{B}$.}

\Proof {(1) If $A$ is closed, then $A=\overline{A}$.\m

(2) If $A=\overline{B}$ for some $B\subset X$, then
$\overline{A}=\overline{\overline{B}}=\overline{B}=A$,
hence $A$ is closed.}

\Remark {$\emptyset$ is open and $X$ is closed.}

\Proof {$X\subset\overline{X}\subset X$, hence $X=\overline{X}$. By Remark 3.6 X is closed.\\
Therefore $\emptyset=X\setminus X$ is open.}

\Remark {$\overline{\emptyset}$ is the smallest closed subset of $X$.}

\Proof {(1) $\overline{\emptyset}$ is closed by Remark 3.6.\m

(2) Let $B$ be a closed subset of $X$. Since $\emptyset\subset B$, we have
$\overline{\emptyset}\subset\overline{B}=B$.}

\Definition {For $x\in X$ we set\m

$\UUU(x):=\gk{U\subset X \mid x\notin\overline{X\setminus U}}$\m

\NI The elements of $\UUU(x)$ are called \K{neighborhoods} of $x$.}

\Remark {Let $A\subset X$ and $x\in X$. Then the following conditions are equivalent:\m

(1) $x\in\overline{A}$.

(2) For every $U\in\UUU(x)$ one has $U\cap A\nleer$.}

\Proof {(1) $\Mpi$ (2): Let $x\in\overline{A}$ and $U\in\UUU(x)$. Then $x\notin\overline{X\setminus U}$.
Assume that $U\cap A\leer$. Then $A\subset X\setminus U$, hence $\overline{A}
\subset\overline{X\setminus U}$. Therefore $x\notin\overline{A}$, a contradiction.\m

(2) $\Mpi$ (1): Assume $x\notin\overline{A}$ and condition (2). From $x\notin\overline{A}
=\overline{X\setminus(X\setminus A)}$ we see that $X\setminus A\in\UUU(x)$.
But $(X\setminus A)\cap A\leer$, a contradiction to (2).}

\Remark {For $x\in X$ the following conditions are equivalent:\m

(1) $x$ is inessential.

(2) $\UUU(x)\leer$.

(3) $X\notin\UUU(x)$.}

\Proof {(1) $\Mpi$ (2): Assume that there exists a neighborhood $U\in\UUU(x)$.
Then $x\notin\overline{X\setminus U}$, hence also $x\notin\overline{\emptyset}$ since
$\overline{\emptyset}\subset\overline{X\setminus U}$. But this means that $x$ is essential.\m

(2) $\Mpi$ (3): Clear.\m

(3) $\Mpi$ (1): Assume $X\notin\UUU(x)$. Then $x\in\overline{X\setminus X}=\overline{\emptyset}$.}

\Remark {A subset $U\subset X$ is open iff $U\in\UUU(x)$ for every $x\in U$.}

\Proof {(1) Let $U$ be open and $x\in U$. Assume that $U\notin\UUU(x)$, i.e. that
$x\in\overline{X\setminus U}$. Since $U$ is open, $\overline{X\setminus U}=X\setminus U$,
hence $x\in X\setminus U$, a contradiction.\m

(2) Assume that $U\in\UUU(x)$ for every $x\in U$ and that $U$ is not open.
Then $\overline{X\setminus U}\setminus(X\setminus U)\nleer$, hence there exists
$x\in\overline{X\setminus U}$ with $x\in U$. By hypothesis $U\in\UUU(x)$, hence
$(X\setminus U)\cap U\nleer$ by Remark 3.10, a contradiction.}

\Remark {Let $x\in X$ and $U\in\UUU(x)$. If $U\subset V\subset X$, then $V\in\UUU(x)$.}

\Proof {By hypothesis, $x\in X\setminus\overline{X\setminus U}
\subset X\setminus\overline{X\setminus V}$, hence $V\in\UUU(x)$.}

\Definition {For $A\subset X$ the \K{interior} of $A$ is defined as\m

$\Moi A:=\gk{x\in X \mid A\in\UUU(x)}$\m

\NI By Remark 3.12 $A$ is open iff $A=\Moi A$.}

\Remark {Let $A\subset X$. Then:\m

(1) $\Moi A=X\setminus\overline{X\setminus A}$.

(2) $\overline{A}=X\setminus\Moi(X\setminus A)$.}

\Proof {(1) $x\in\Moi A\gd A\in\UUU(x) \gd x\notin\overline{X\setminus A}$.\m

(2) From (1), substituting $X\setminus A$ for $A$, we have\m

$\Moi (X\setminus A)=X\setminus\overline{X\setminus(X\setminus A)}=
X\setminus\overline{A}$\m

\NI hence $\overline{A}=X\setminus\Moi(X\setminus A)$.}

\Proposition {For $x\in X$ and $U\subset X$ the following conditions are equivalent:\m

(1) $U\in\UUU(x)$.

(2) There exists an open set $W$ such that $x\in W\subset U$.}

\Proof {(1) $\Mpi$ (2): Set $W:=X\setminus\overline{X\setminus U}$. Then $W$ is open
by Remark 3.6 and from $U\in\UUU(x)$ we have $x\in X\setminus\overline{X\setminus U}=W
\subset X\setminus(X\setminus U)=U$.\m

(2) $\Mpi$ (1): Clear from Remarks 3.12 and 3.13.}

\Remark {The following conditions are equivalent:\m

(1) $X$ is a topological space.

(2) $\UUU(x)$ is a filter on $X$ for every $x\in X$.}

\Proof {Clear. Notice that (2) implies that $X\in\UUU(x)$ for every $x\in X$, hence all
points of $X$ are essential. Cf. Proposition 3.38.}

\Remark {Let $x\in X$. Then every neighborhood of $x$ contains a minimal neighborhood of $x$.}

\Definition {For $x\in X$ let $\MMM(x):=\Min\UUU(x)$ be the set of all minimal
neighborhoods of $x$.

$x$ is inessential iff $\MMM(x)\leer$.}

\Definition {For $A\subset X$ let $\MMM(A):=\Movg_{a\in A}\MMM(a)$.

In particular $\MMM(X)=\Movg_{x\in X}\MMM(x)$.}

\Remark {For $x\in X$ one has\m

$\UUU(x)=\gk{U\subset X \mid \text{ there exists } M\in\MMM(x)
\text { such that } M\subset U}$}

\Remark {Let $x\in X$. Then every element of $\MMM(x)$ is open.}

\Proof {This follows from Proposition 3.16.}

\Definition {Let $f:X\Mp Y$ be a mapping, $x\in X$ and $y:=f(x)$.

$f$ is \K{continuous} in $x$ if for every $M\in\MMM(x)$ there exists $N\in\MMM(y)$
such that $f(M)\subset N$.

$f$ is \K{continuous} if it is continuous in every point of $X$.}

\Remark {Let $f:X\Mp Y$ be a mapping. Then $f$ is continuous in every inessential point of $X$.}

\Definition {Let $x\in X$ and $A\subset X$. We say that $A$ \K{converges} to $x$
and write $A\Mp x$, if there exists $M\in\MMM(x)$ such that $A\subset M$.

We set $\CCC(x):=\gk{A\subset X \mid A\Mp x} = \Movg_{M\in\MMM(x)}\PPP(M)$.}

\Remark {Let $x\in X$. Then $\CCC(x)$ has the following properties:\m

(1) $A\subset B\in\CCC(x) \Mpi A\in\CCC(x)$.\m

(2) $A\Mp x\Mpi A\cup x\Mp x$.}

\Proof {(1) Assume $A\subset B\in\CCC(x)$. Then there exists $M\in\MMM(x)$ such that
$B\subset M$. Hence also $A\subset M$, therefore $A\in\CCC(x)$.\m

(2) Let $A\Mp x$. Then there exists $M\in\MMM(x)$ such that $A\subset M$.

But $M\in\UUU(x)$, therefore $x\in M$, so that $A\cup x\subset M$.
Thus $A\cup x\Mp x$.}

\Remark {For $x\in X$ the following conditions are equivalent:\m

(1) $x$ is essential.\m

(2) $\emptyset\Mp x$.\m

(3) $x\Mp x$.\m

(4) $\CCC(x)\nleer$.}

\Proof {(1) $\Mpi$ (2): Since $x$ is essential, there exists $A\subset X$ such that
$A\Mp x$. Since $\emptyset\subset A$, this implies $\emptyset \Mp x$.\m

(2) $\Mpi$ (3): If $\emptyset\Mp x$, then by Remark 3.26 also $\emptyset\cup x=x\Mp x$.\m

(3) $\Mpi$ (4) $\Mpi$ (1): Clear.}

\Corollarynk {From Remarks 3.26 and 3.27 one sees that, if $x$ is an\\essential point, then
$\CCC(x)$ is an \K{abstract simplicial complex} on $X$\\
(cfr. Kozlov [10, p. 7] and Barmak [9, p. 151]).}

\Proposition {Let $f:X\Mp Y$ be a mapping and $x\in X$. Then the\\following statements
are equivalent:\m

(1) $f$ is continuous in $x$.\m

(2) $A\Mp x\Mpi f(A)\Mp f(x)$.}

\Proof {Let $y:=f(x)$.\m

(1) $\Mpi$ (2): Assume that $f$ is continuous in $x$ and that $A\Mp x$. Then there exists
$M\in\MMM(x)$ such that $A\subset M$, and by the continuity of $f$ in $x$ there exists
$N\in\MMM(y)$ such that $f(M)\subset N$. Then also $f(A)\subset N$ and this implies
that $f(A)\Mp y$.\m

(2) $\Mpi$ (1): Take $M\in\MMM(x)$. Then $M\Mp x$, hence, by hypothesis (2),
$f(M)\Mp y$. Therefore there exists $N\in\MMM(y)$ such that $f(M)\subset N$.}

\Proposition {Let $f:X\Mp Y$ and $g:Y\Mp Z$ be mappings and let $x\in X$.
Assume that $f$ is continuous in $x$ and that $g$ is continuous in $f(x)$.

Then $g\circ f$ is continuous in $x$.}

\Proof {Let $A\Mp x$. Then $f(A)\Mp f(x)$ since $f$ is continuous in $x$, and
$g(f(A))\Mp g(f(x)$ since $g$ is continuous in $f(x)$.}

\Remark {Let $T$ be a finite \K{topological} space and $t\in T$. Then:\m

(1) $\MMM(t)=\gk{U_t}$.\m

(2) $\CCC(t)=\PPP(U_t)$.

\zweim Hence $Q\Mp t$ iff $Q\subset U_t$.}

\Proof {Clear.}

\Lemma {Let $T$ be a finite topological space and $f:T\Mp X$ a mapping. Then for $t\in T$
and $x:=f(t)$ the following conditions are equivalent:\m

(1) $f$ is continuous in $t$.

(2) There exists $M\in\MMM(x)$ such that $f(U_t)\subset M$.

(3) $f(U_t) \Mp x$.}

\Proof {(1) $\Mpi$ (2): Let $f$ be continuous in $t$. Since $U_t\Mp t$, we have
$f(U_t)\Mp x$. This means that there exists $M\in\MMM(x)$ such that $f(U_t)\subset M$.\m

(2) $\Mpi$ (3): By definition.\m

(3) $\Mpi$ (1): Let $A\Mp t$. Then $A\subset U_t$, hence $f(A)\subset f(U_t)\Mp x$,
therefore $f(A)\Mp x$.}

\Definition {A mapping $f:X\Mp Y$ is said to be \K{open}, if for every open subset $U\subset X$
its image $f(U)$ is open in $X$.}

\Lemma {Let $f:X\Mp Y$ be a mapping. The following conditions are equivalent:\m

(1) $f$ is open.

(2) For every $x\in X$ and every $M\in\MMM(x)$ the image $f(M)$ is open in $Y$.}

\Proof {(1) $\Mpi$ (2): Clear, since every $M\in\MMM(x)$ is open by Remark 3.22.\m

(2) $\Mpi$ (1): Let $U\subset X$ be open. Take $y\in f(U)$, i.e. $y=f(x)$ for some $x\in U$.
By hypothesis $U\in\UUU(x)$ and by Remark 3.21 there exists $M\in\MMM(x)$ such that
$M\subset U$.

By (2) then $f(M)$ is open in $Y$. Since $x\in M$, we have $y\in f(M)$, hence
$f(M)\in\UUU(y)$ and therefore, since $f(M)\subset f(U)$, also $f(U)\in\UUU(y)$.}

\Corollary {Let $f:X\Mp Y$ be a mapping. The following\\conditions are equivalent:\m

(1) $f$ is continuous and open.

(2) For every $x\in X$ and every $M\in\MMM(x)$ one has $f(M)\in\MMM(f(x))$.}

\Proof {(1) $\Mpi$ (2): Let $x\in X$ and $y:=f(x)$. Since $f$ is open, from Lemma 3.34
we have $f(M)\in\UUU(y)$. This implies that there exists $K\in\MMM(y)$ with
$K\subset f(M)$. But $f$ is also continuous, therefore there exists $N\in\MMM(y)$
such that $f(M)\subset N$. Then $K\subset f(M)\subset N$, thus $K=N$ by
minimality, hence also $f(M)=N\in\MMM(y)$.\m

(2) $\Mpi$ (1): By Lemma 3.34 $f$ is open. $f$ is cleary continuous.}

\Proposition {Let $x\in X$, $M\in\MMM(x)$ and $y\in M$.

Then $M\in\MMM(y)$ or $M\setminus x\in\UUU(y)$.}

\Proof {$M$ is open by Remark 3.22, therefore $M\in\UUU(y)$.\\Assume that $M\notin\MMM(y)$.

Then there exists $N\in\MMM(y)$ with $N\subset M$. Assume that $x\in N$.

But also $N$ is open by Remark 3.22, hence $N\in\UUU(x)$. Now $M\in\MMM(x)
=\Min\UUU(x)$, therefore $M=N\in\MMM(y)$, a contradiction since we assumed
that $M\notin\MMM(y)$.

Therefore $N\subset M\setminus x$ and this implies $M\setminus x\in\UUU(y)$.}

\Definition {A point $x\in X$ is said to be \K{regular}, if $|\MMM(x)|=1$.\\
A non-regular point is called \K{singular}.

Notice that a regular point is necessarily essential.}

\Proposition {Let $x$ be an essential point of $X$.

Then the following conditions are equivalent:\m

(1) $x$ is regular.

(2) $\UUU(x)$ is a filter on $X$.

(3) $U,V\in\UUU(x) \Mpi U\cap V\in\UUU(x)$.

(4) $A,B\Mp x\Mpi A\cup B\Mp x$.}

\Proof {(1) $\Mpi$ (2): Let $\MMM(x)=\gk{M}$. Then $\UUU(x)=\widehat{M}$ and
this is a filter (since $M\nleer$).\m

(2) $\Mpi$ (1): Let $M:=\Modg_{U\in\UUU(x)}U$. Then, since by hypothesis $\UUU(x)$
is a filter and finite, $M\in\UUU(x)$, hence $\MMM(x)=\Min\UUU(x)=\gk{M}$.\m

(2) $\gd$ (3): Clear (since $x$ is essential).\m

(1) $\Mpi$ (4): Assume $\MMM(x)=\gk{M}$ and let $A,B\Mp x$.

Then necessarily $A,B\subset M$, hence also $A\cup B\subset M$, thus $A\cup B\Mp x$.\m

(4) $\Mpi$ (1): Let $M,N\in\MMM(x)$. By hypothesis $M\cup N\Mp x$ and the maximality
of $M$ and $N$ implies that $M=M\cup N$ and $N=M\cup N$, hance $M=N$.}

\Corollary {The following conditions are equivalent:\m

(1) $X$ is a topological space.

(2) $X$ does not contain inessential points and for every $x\in X$\\\zweim one has
$A,B\Mp x \Mpi A\cup B\Mp x$.}

\Definition {A mapping $f:X\Mp Y$ is said to be \K{combinatorially\\[0.6ex] continuous}
in $x\in X$, if for every $V\in\UUU(f(x))$ one has $f\hme(V)\in\UUU(x)$.

$f$ is called \K{combinatorially continuous}, if it is continuous in every\\point of $X$.}

\Remark {Let $f:X\Mp Y$ be a mapping and $x\in X$.\\Then $f$ is combinatorially
continuous in $x$ iff for every $V\in\UUU(f(x))$\\there exists $U\in\UUU(x)$ with $f(U)\subset V$.}

\Proof {This follows from $f(U)\subset V\gd U\subset f\hme(V)$ and Remark 3.13:\m

(1) Assume that $f$ is combinatorially continuous in $x$ and let $V\in\UUU(f(x))$.
By hypothesis one has $U:=f\hme(V)\in\UUU(x)$. Then $f(U)=f(f\hme(V))\subset V$.\m

(2) Let the condition (2) be true. Take $V\in\UUU(f(x))$. By hypothesis there exists $U\in\UUU(x)$
such that $f(U)\subset V$. Then $U\subset f\hme(f(U))\subset f\hme(V)$, hence
$f\hme(V)\in\UUU(x)$.}

\Proposition {Let $f:X\Mpi Y$ be a mapping. Then the following conditions are equivalent:\m

(1) $f$ is combinatorially continuous.

(2) For every $x\in X$ and every $V\in\UUU(f(x))$ there exists $U\in\UUU(x)$\\
\zweim such that $f(U)\subset V$.

(3) For every open subset $V$ of $Y$ the preimage $f\hme(V)$ is open in $X$.

(4) For every closed subset $B$ of $Y$ the preimage $f\hme(B)$ is closed in $X$.

(5) For every $A\subset X$ one has $f(\overline{A})\subset \overline{f(A)}$.}

\Proof {(1) $\gd$ (2): Remark 3.41.\m

(2) $\Mpi$ (3): Let $V$ be open in $Y$ and $x\in f\hme(V)$. Then $f(x)\in V$,\\hence
$V\in\UUU(f(x))$. By (2) there exists $U\in\UUU(x)$ such that $f(U)\subset V$,
i.e. $U\subset f\hme(V)$. Therefore $f\hme(V)\in\UUU(x)$.\m

(3) $\gd$ (4): This follows from $f\hme(Y\setminus B)=X\setminus f\hme(B)$.\m

(4) $\Mpi$ (5): Let $x\in\overline{A}$ and set $C:=f\hme(\overline{f(A)})$. Then $A\subset C$
and, by (4), $C$ is a closed subset of $X$. Therefore $\overline{A}\subset
\overline{C}=C$, hence $f(\overline{A})\subset f(C)=f(f\hme(\overline{f(A)}))\subset \overline{f(A)}$.

(5) $\Mpi$ (1): Let $V\in\UUU(f(x))$ and suppose that $f\hme(V)\notin\UUU(x)$.
This means that $x\in\overline{X\setminus f\hme(V)}$, hence
$f(x)\in f(\overline{X\setminus f\hme(V)})\subset \overline{f(X\setminus f\hme(V))}$,\\
thus $f(X\setminus f\hme(V))\cap V\nleer$.

Therefore there exists $a\in V$ such that $a=f(b)$ for some $b\in X\setminus f\hme(V)$,
i.e. $a=f(b)\notin V$, a contradiction.}

\Remark {Condition (5) in Prop. 3.42 is the defining property commonly used
in combinatorics for mappings between closure spaces.\\Cf. Ern\'e [3, p. 174-175].}
\kapitel{4}
\Ab{4. The topological resolution}

\Standing {Let $X$ and $Y$ be \K{finite} closure spaces.}

\Definition {$\Top X:=\gk{xM \mid x\in X \text{ and } M\in\MMM(x)}$.

Recall that here $xM$ is a short-cut for the ordered pair $(x,M)$.\m

We define a quasiorder (hence a topology) on $\Top X$ by\m

$xM\le yN :\gd N\subset M$\m

\NI We have a natural projection $\Mfa {\pi:\Top X} {X}{xM} {x}$\m

\NI By Definition 3.19 the image of $\pi$ coincides with the set of all essential points of $X$.

Therefore $\pi$ is surjective iff every point of $X$ is essential, i.e. iff \; $\overline{\emptyset}\leer$.

We call the topological space $\Top X$ the \K{topological resolution} of the\\closure space $X$.}

\Remark {For $xM\in\Top X$ one has

$U_{xM} = \gk{yN \in\Top X\mid N\subset M}$}

\Proof {For $y\in X$ and $N\in\MMM(y)$ one has (by Proposition 2.1):\m

$yN\in U_{xM} \gd yN\ge xM\gd N\subset M$}

\Definition {By Remark 4.3 the neighborhood $U_{xM}$ depends only on $M$, not on $x$,
in the sense that if $M\in\MMM(x)\cap\MMM(y)$, then $U_{xM}=U_{yM}$.\\We introduce
the following notation:\m

For $A\subset X$ we set\quad
$[A]:=\gk{yN\in\Top X \mid N\subset A}$.\m

\NI For $x\in X$ and $M\in\MMM(x)$ then $U_{xM}=[M]$,
hence\; $\UUU(xM)=\widehat{U_{xM}}=\widehat{[M]}$.}

\Theorem {For $A\subset X$ we have $\Moi A=\pi([A])$.}

\Proof {(1) Let $x\in\Moi A$, i.e. $A\in\UUU(x)$. Then there exists $M\in\MMM(x)$ with
$M\subset A$, thus $xM\in [A]$, therefore $x=\pi(xM)\in\pi([A])$.\m

(2) Let $x\in\pi([A])$. Then there exists $yM\in [A]$ such that $x=\pi(yM)=y$.
From $M\subset A$ it follows that $A\in\UUU(x)$, hence $x\in\Moi A$.}

\Proposition {Let $xM\in\Top X$. Then $\pi(U_{xM})=\pi([M])=M$.}

\Proof {This follows from Theorem 4.5, since $M=\Moi M$ by Remark 3.22.}

\Theorem {The natural projection $\pi:\Top X\Mp X$ is continuous\\and open.}

\Proof {(1) Let $xM\in\Top X$. Then $M\in\MMM(x)$ and $\pi(U_{xM})=M$ by Proposition 4.6.
By Lemma 3.32 $\pi$ is continuous.\m

(2) $\pi$ is open by Lemma 3.34, Proposition 4.6 and Remark 3.22.}

\Proposition {For a mapping $f:X\Mp Y$ consider the composition\\
$\Top X\Mpb{\pi} X\Mpb{f} Y$.

Then $f$ is continuous iff $f\circ\pi$ is continuous.}

\Proof {(1) If $f$ is continuous, then $f\circ \pi$ is continuous by Proposition 3.30.\m

(2) Assume that $f\circ\pi$ is continuous. Let $x\in X$ and $M\in\MMM(x)$.

Then $xM\in\Top X$ and by Lemma 3.32 (and the continuity of $f\circ\pi$) there exists
$N\in\MMM(f(\pi(xM)))=\MMM(f(x))$ such that $(f\circ\pi)(U_{xM})\subset N$.

But $(f\circ\pi)(U_{xM})=f(M)$ by Proposition 4.6, hence $f(M)\subset N$.
This means that $f$ is continuous in $x$.}

\Remark {If $X$ is a topological space, then the natural projection\\$\pi:\Top X\Mp X$ is a
homeomorphism.}

\Proof {Immediate. Notice that in this case $\Top X=\gk{xU_x \mid x\in X}$.}

\Lemma {Let $A,B\subset X$. Then:\m

(1) $A\subset B \Mpi [A]\subset [B]$.

(2) $[A\cap B]=[A]\cap [B]$.

(3) $[X]=\Top X$.}

\Proof {(1) $xM\in [A]\Mpi M\subset A\Mpi M\subset B\Mpi xM\in [B]$.\m

(2) $xM\in [A\cap B]\gd M\subset A\cap B\gd M\subset A \text{ and } M\subset B\\[1ex]
\zweim\gd xM\in [A]\cap [B]$.\m

(3) Let $xM\in\Top X$. Then $M\subset X$, hence $xM\in [X]$.}

\Remark {(1) $[A]$ is open in $\Top X$ for every $A\subset X$.\m

(2) The neighborhood filter $\UUU(xM)$ is the set of all $O\subset\Top X$ with the property
that there exists $A\subset X$ such that $xM\in [A]\subset O$.\m

(3) The families $\gk{[M] \mid M\in\MMM(X)}$ and $\gk{[A] \mid A\subset X}$ constitute
both a basis for the open subsets of $\Top X$.}

\Proof {(1) Let $xM\in [A]$. Then $M\subset A$, hence $xM\in [M]\subset [A]$.

Since $[M]\in\UUU(xM)$, this implies $[A]\in\UUU(xM)$.\m

(2) Let $O\in\UUU(xM)=\widehat{U_{xM}}=\widehat{[M]}$. Then $xM\in [M]\subset O$.\m

If viceversa $xM\in [A]\subset O$, then by (1) $O\in\UUU(xM)$.\m

(3) Follows from (1) and (2).}

\Proposition {Let $W\subset X$. Then $W$ is open in $X$ iff there exists $A\subset X$
such that $W=\pi([A])$.}

\Proof {(1) If $W$ is open in $X$, then $W=\Moi W=\pi([W])$.\m

(2) If $W=\pi([A])$, then $W$ is open, since $[A]$ is open in $\Top X$
and $\pi$ is an open mapping.}

\Corollary {A subset of $X$ is open iff it is the image under $\pi$ of an open subset of $\Top X$.}

\Remark {(1) Let $X$ be a set. It is well known (see e.g. Ern\'e [2] or Ihringer [7, p.36]) that every
family $\EEE$ of subsets of $X$ with the property that $A,B\in\EEE$ implies $A\cup B\in\EEE$
can be considered as the family of open subsets of a closure space.

From this follows that, if $T$ is a finite topological space, $X$ a set and\\$p:T\Mp X$ a mapping,
then we can define a closure space structure on $X$ using the family
$\EEE:=\gk{p(O) \mid O \text { open in } T}$. The inessential points are
the elements of $X\setminus p(T)$.\m

(2) From the foregoing discussion, in particular from Corollary 4.13,\\we conclude that every
finite closure space may be obtained in this way.}

\Lemma {Let $f:X\Mp Y$ be an isomorphism (i.e. a bijective mapping which is continuous
in both directions), $x\in X$ and $M\in\MMM(x)$. Then $f(M)\in\MMM(f(x))$.}

\Proof {Let $g:=f\hme$ and $y:=f(x)$.

$f$ is continuous, therefore there exists $N\in\MMM(y)$ such that $f(M)\subset N$.

$g$ is continuous, therefore there exists $K\in\MMM(x)$ such that $g(N)\subset K$.

Then $M=g(f(M))\subset g(N)\subset K$, and this implies $K=M$ by minimality of $M$.
Therefore $g(N)=M$ and $f(M)=N\in\MMM(y)$.}

\Proposition {Let $f:X\Mp Y$ be an isomorphism. Then the mapping\m

$\Mfa {F:\Top X}{\Top Y}{xM}{f(x)f(M)}$\m

\NI is well defined and a homeomorphism.}

\Proof {From Lemma 4.15 it follows that $F$ is well defined and bijective.
It is also clear that the mappings $F$ and $F\hme$ are monotone and therefore continuous
(cf. Proposition 2.1).}

\Remark {Let $x\in X$. Then $\MMM(x)=\gk{\pi(U_t) \mid t\in\pi\hme(x)}$.}

\Proof {(1) Let $M\in\MMM(x)$. Then $t:=xM\in\pi\hme(x)$ and $M=\pi(U_t)$\\
by Proposition 4.6 and Definition 4.4.\m

(2) Let $t\in\pi\hme(x)$. Then $t=xM$ for some $M\in\MMM(x)$.\\Therefore
$\pi(U_t)=M\in\MMM(x)$ by Proposition 4.6.}

\Corollary {Let $x\in X$ and $A\subset X$. Then the following conditions are equivalent:\m

(1) $A\Mp x$.

(2) There exists $t\in\pi\hme(x)$ with $A\subset\pi(U_t)$.}

\Remark {We are now able to characterize continuous resp. combinatorially continuous mappings
between finite closure spaces in terms of the topological resolutions.}

\Proposition {Let $f:X\Mp Y$ be a mapping, $x\in X$ and $y:=f(x)$.

Then the following conditions are equivalent:\m

(1) $f$ is continuous in $x$.

(2) For every $t\in\pi\hme(x)$ there exists $s\in\pi\hme(y)$ such that
$f(\pi(U_t))\subset\pi(U_s)$.}

\Proposition {A mapping $f:X\Mp Y$ is combinatorially continuous iff for every open
subset $P$ of $\Top Y$ there exists an open subset $O$ of $\Top X$ such that
$f\hme(\pi(P))=\pi(O)$.}

\Proof {This follows from Propositions 3.42 and 4.12.}
\kapitel{5}
\Ab{5. Regular mappings}

\Standing {Let $X$ and $Y$ be \K{finite} closure spaces and\\$f:X\Mp Y$ a mapping.
$x\in X$ where not otherwise indicated.}

\Definition {Let $y:=f(x)$. $f$ is said to be \K{regular} in $x$, if for every\\
$M\in\MMM(x)$ there exists $K\in\MMM(y)$ such that $\widehat{f(M)}\cap \UUU(y)=\widehat{K}$.

In this case $K$ is uniquely determined since $\widehat{K}=\widehat{L}$ implies $K=L$.

$f$ is called \K{regular} (on $X$) if it is regular in every point of $X$}

\Lemma {Let $A\subset X$ be such that $\widehat{A}\cap\UUU(x)=\widehat{K}$
for some $K\in\MMM(x)$. Then $A\subset K$ and there are no other elements of $\MMM(x)$
containing $A$.}

\Proof {(1) We have in particular $K\in\widehat{A}\cap \UUU(x)$, hence $A\subset K$.\m

(2) Let $L\in\MMM(x)$ be such that $A\subset L$. Then $L\in\widehat{A}\cap\UUU(x)=\widehat{K}$,
hence $L\supset K$. Since $L,K\in\MMM(x)$, this implies $L=K$.}

\Remark {Let $f$ be regular in $x$. Then $f$ is continuous in $x$.}

\Proof {Let $M\in\MMM(x)$ and $y:=f(x)$. By hypothesis there exists $K\in\MMM(y)$
such that $\widehat{f(M)}\cap\UUU(y)=\widehat{K}$.

Then $f(M)\subset K$ by Lemma 5.3.}

\Proposition {Let $f(x)$ be a regular point of $Y$ and assume that $f$ is\\continuous in $x$.
Then $f$ is regular in $x$.}

\Proof {Since $y:=f(x)$ is a regular point, we have $\UUU(y)=\widehat{K}$ for the unique element $K$
of $\MMM(y)$. Let $M\in\MMM(x)$. Since $f$ is continuous, we have\\$f(M)\subset K$.

But then $\widehat{K}\subset\widehat{f(M)}$, hence $\widehat{f(M)}\cap\UUU(y)
=\widehat{f(M)}\cap\widehat{K}=\widehat{K}$.}

\Corollary {Let $Y$ be a topological space.

Then $f$ is regular iff $f$ is continuous.}

\Remark {Let $M\in\MMM(x)$. Then $\widehat{M}\subset\UUU(x)$, hence
$\widehat{M}\cap\UUU(x)=\widehat{M}$.}

\Proposition {Let $f$ be continuous and open. Then $f$ is regular.}

\Proof {Let $M\in\MMM(x)$ and $y:=f(x)$. By Corollary 3.35\\$K:=f(M)\in\MMM(y)\subset\UUU(y)$.

Therefore by Remark 5.7 we have $\widehat{f(M)}\cap\UUU(y)=\widehat{f(M)}=\widehat{K}$.}

\Corollary {The natural projection $\pi:\Top X\Mp X$ is regular.}

\Remark {Let $F:\Top X\Mp \Top Y$ be a continuous mapping such that the diagram

\ein{5} $\xymatrix
{\Top X\ar[r]^F\ar[d]_\pi & \Top Y\ar[d]^\pi\\
X\ar[r]^f & Y}$\sm

\NI commutes. Then $f$ is continuous.}

\Proof {This follows from Proposition 4.8.}

\Theorem {If $f$ is regular, then there exists a unique continuous\\mapping $F:\Top X\Mp\Top Y$
such that the diagram

\ein{5} $\xymatrix
{\Top X\ar[r]^F\ar[d]_\pi & \Top Y\ar[d]^\pi\\
X\ar[r]^f & Y}$\sm

\NI commutes. $F$ is defined in the following way:\m

Let $xM\in\Top X$ and $y:=f(x)$. By regularity of $f$ there exists $K\in\MMM(y)$ such that
$\widehat{f(M)}\cap\UUU(y)=\widehat{K}$.

Then we set $F(xM):=yK$.}

\Proof {Let $x,M,y,K$ be as in the statement of the theorem. Evidently $F$ is well defined
as a mapping.\m

(1) Commutativity of the diagram is immediate:

$f(\pi(xM))=f(x)=y$

$\pi(F(xM))=\pi(yK)=y$\m

(2) We show the continuity of $F$ in $x$. As in the proof of Remark 5.4 (or in Lemma 5.3) we have
$f(M)\subset K$, hence $M\subset f\hme(K)$. Therefore $xM\in[f\hme(K)]$ and
$[f\hme(K)]$ is an  open neighborhood of $xM$ in $\Top X$. Since $[K]$ is the minimal
neighborhood of $yK=F(xM)$ in $\Top Y$, it suffices to show that
$F([f\hme(K)])\subset [K]$.

Let $sN\in[f\hme(K)]$. Then $N\subset f\hme(K)$, i.e. $f(N)\subset K$.
By regularity of $f$ we have $\widehat{f(N)}\cap\UUU(f(s))=\widehat{L}$
for some $L\in\MMM(f(s))$.

Now $N\in\MMM(s)$, hence $s\in N\subset f\hme(K)$, so that $f(s)\in K$, therefore\\
$K\in\UUU(f(s))$. Since $f(N)\subset K$, we have $K\in\widehat{f(N)}$, hence
$K\in\widehat{L}$, i.e. $L\subset K$. This implies $F(sN)=f(s)L\in[K]$.\m

(3) Unicity of $F$: Let $G:\Top X\Mp\Top Y$ be continuous and such that the diagram

\ein{5} $\xymatrix
{\Top X\ar[r]^G\ar[d]_\pi & \Top Y\ar[d]^\pi\\
X\ar[r]^f & Y}$\sm

\NI is commutative. Let $xM\in\Top X$ and $G(xM)=yR$. Since $G$ is continuous, we must
have $G(U_{xM})\subset U_{yR}$, i.e. (by Definition 4.4) $G([M])\subset [R]$.

Using Proposition 4.6 we have\m

$f(M)=f(\pi([M]))=\pi(G([M]))\subset\pi([R])=R$\m

\NI On the other hand, because $f$ is regular, we have $f(M)\subset K$.
But then $R=K$ by Lemma 5.3, hence $G(xM)=F(xM)$.}
\newpage\kapitel{5}
\Ab{References}\m

\centerline{\K{Closure spaces}}\mm

\begin{biblio}
\aitem{2}{3084} {G. Birkhoff} Lattice theory. AMS 1967.

\aitem {3} {16134} {M. Ern\'e} Einführung in die Ordnungstheorie. Bibl. Inst. 1982.

\aitem {4} {24541} {M. Ern\'e} Closure. Contemporary Mathematics 486 (2009), 163-238.

\aitem {5} {7331} {J. Eschgfäller} Almost topological spaces.\\
Ann. Univ. Ferrara 30 (1984), 163-183.

\aitem {7} {26139} {B. Ganter} Diskrete Mathematik - geordnete Mengen. Springer 2013.

\aitem {8} {14111} {B. Ganter/R. Wille} Formale Begriffsanalyse. Springer 1996.

\aitem {9} {1181} {T. Ihringer} Allgemeine Algebra. Teubner 1988.

\aitem {11} {21603} {G. Nöbeling} Grundlagen der analytischen Topologie. Springer 1954.

\end{biblio}

\bigskip

\centerline{\K{Finite topological spaces and topological combinatorics}}

\begin{biblio}
\aitem {1} {29028} {J. Barmak} Algebraic topology of finite topological spaces and\\applications.
Springer 2011.

\aitem {2} {16061} {A. Brini} Combinatoria e topologia.\\
Boll. UMI Mat. Soc. Cultura Dicembre 2003, 531-563.

\aitem {10} {29026} {D. Kozlov} Combinatorial algebraic topology.\\Springer 2008.

\aitem {2} {29455} {M. de Longueville} A course in topological combinatorics.\\Springer 2013.

\aitem {9} {29456} {J. Matou\v{s}ek} Using the Borsuk-Ulam theorem. Springer 2003.

\aitem {11} {29040} {J. May} Finite topological spaces. Internet 2008, 13p.

\aitem {12}{7193} {R. Stong} Finite topological spaces. Trans. AMS 123 (1966), 325-340.
\end{biblio}

\end{document}